\documentclass[12pt]{article}

\title{Removability of time-dependent singularities \\
in the heat equation}

\author{ 
Jin Takahashi
	\footnote{
	Department of Mathematics, 
	Tokyo Institute of Technology, 
	Meguro-ku, Tokyo 152-8551, Japan}
	\footnote{
	E-mail: takahashi.j.ab@m.titech.ac.jp
	}
\and
Eiji Yanagida 
	\footnotemark[1]
}

\usepackage{amsmath,amsthm,amssymb,comment}

\makeatletter

\newtheorem{thm}{Theorem}[section]
\newtheorem{lem}{Lemma}[section]

\newtheorem{pro}{Proposition}[section]

\newcommand{\Sect}[1]{\setcounter{equation}{0}\section{#1}\quad}

\@addtoreset{equation}{section}

\makeatother

\newcommand{\Qed}{\quad\hbox{\rule[-2pt]{3pt}{6pt}}\par\bigskip}
\newcommand{\R}{\mathbb{R}}
\newcommand{\eps}{\varepsilon}

\newcommand{\loc}{{\rm loc}}
\newcommand{\rank}{\mathop{\rm rank}}
\newcommand{\tu}{\tilde u}
\newcommand{\txi}{\tilde \xi}
\newcommand{\tXi}{\widetilde \Xi}
\newcommand{\bs}{\mathbf{s}}
\newcommand{\sm}{s_1,s_2,\ldots,s_m}
\theoremstyle{definition}

\begin{document}
\date{}

\maketitle

\begin{abstract}
We consider solutions of the linear heat equation with
time-dependent singularities.
It is shown that
if a singularity is weaker than the order of 
the fundamental solution
of the Laplace equation, 
then it is removable.  We also consider the removability of
higher dimensional singular sets.  
An example of a non-removable singularity is given, which
implies the optimality of the condition for removability.
\end{abstract}

\Sect{Introduction}\label{sect:intro}

Removability of singularities of solutions is an interesting and
important problem in partial differential equations.
For the Laplace equation, the removability of a singular point is
defined as follows. Let $u$ be a solution of 
\[
	\Delta u=0 \quad \mbox{ in } \Omega \setminus \{\xi_0 \}, 
\]
where $\Omega$ is a domain in $\R^N$ and $\xi_0 \in \Omega$. 
We say that $\xi_0$ is a removable 
singularity if there exists a classical solution $\tu$ of the
Laplace equation in $\Omega$ such that
\[
	\tu \equiv u \quad \mbox{ in } \Omega \setminus \{\xi_0 \}.
\]
It is well known \cite{F} that for $N\geq 3$, 
the singular point $\xi_0$ is 
removable if and only if
\[
	|u(x)|=o(|x-\xi_0|^{2-N}) \quad \mbox{as} \quad x \to \xi_0.
\]
For nonlinear elliptic equations, the removability of a singularity has
been studied in many papers 
and various interesting results have been obtained
(see, e.g.,  Brezis-Veron~\cite{BV},  Gidas-Spruck~\cite{GS}, Veron~ \cite{V}).

Similarly, for the heat equation
\[
	u_t=\Delta u \quad 
		\mbox{ in } \Omega \setminus \{\xi_0 \} \times (0,T)
\]
with $N\geq 3$ and $T>0$, 
Hsu \cite{Hsu} proved recently that the singular point $\xi_0$ is 
removable if and only if 
\[
	|u(x,t)|=o(|x-\xi_0|^{2-N}) 
	\quad \mbox{as} \quad x\to \xi_0
\]
for every $t\in(0,T)$. 
Later, Hui \cite{Hui} gave a simpler proof for this result.
In \cite{Hirata}, Hirata extended Hsu and Hui's result
to a semilinear parabolic equation
of the form 
\[
u_t=\Delta u+|u|^{p-1}u
\]
 with $p<N/(N-2)$.
See also Sato-Yanagida~\cite{SY} for non-removable 
 singularities of this equation.


In this paper, we consider the case where a singular point may move in time
and study its removability for the heat equation.
More precisely, we formulate our problem as follows.
For $T>0$ fixed, 
let $\xi:[0,T] \rightarrow \R^N$ be a continuous function, 
and $\Gamma \subset \R^{N+1}$  be a curve given by
\[
	\Gamma 
	:= \{ (x,t) \in \R^{N+1} \, : \, x=\xi(t), \; t \in (0,T) \}.
\]
We take a domain $\Omega \subset \R^{N}$ such that 
 $\xi(t) \in \Omega$ for $t \in [0,T]$, 
and define 
\[
	D 
	:= \{
		 (x,t) \in \R^{N+1} \,:\,
		 x\in \Omega \setminus\{\xi(t)\}, \; t\in(0,T) 
		\}.
\]
For a solution of 
\begin{equation} \label{eq:main}
	u_t=\Delta u \quad \mbox{ in } D, 
\end{equation}
the singularity at $x=\xi(t)$ is said to be removable
if there exists
a function $\tu$ which satisfies the heat equation 
in $\Omega \times(0,T)$ in the classical sense and 
$\tu \equiv u$ on $D$.

Our first result gives a
condition for the removability of such a (moving) singularity.

\begin{thm}\label{th:remove0} \ 
Let $N\geq3$. Suppose that $\xi$ is H\"older continuous
with exponent $1/2$ and that $u$ satisfies $(\ref{eq:main})$ 
in the classical sense.
Then the singularity of $u$ at $x=\xi(t)$ is removable if and only if
for any $0<t_1<t_2<T$ and $0<\eps <1$ 
there exists $0<r<1$ depending on $t_1,t_2,\eps$ such that
\begin{equation}\label{eq:eps0}
	|u(x,t)| 
	\leq 
	\frac{ \eps }{ |x-\xi(t)|^{N-2} },
		\quad 0<|x-\xi(t)|<r
\end{equation}
for any $t \in [t_1 , t_2]$.
\end{thm}

\begin{thm}\label{th:2remove0} \ 
Let $N=2$. 
Suppose that $\xi$ is H\"older continuous
with exponent $1/2$ 
and that $u$ satisfies $(\ref{eq:main})$ 
in the classical sense. 
Then the singularity of $u$ at $x=\xi(t)$ is removable if and only if 
for any $0<t_1<t_2<T$ and $0<\eps <1$ the function $u$ satisfies
\begin{equation}\label{eq:2eps0}
	|u(x,t)| 
	\leq 
	\eps \log\frac{1}{|x-\xi(t)|},
		\quad 0<|x-\xi(t)|<\eps
\end{equation}
for any $t \in [t_1 , t_2]$.
\end{thm}

Here we note that
for $N=1$, if we define $\tu$ by 
\[
	\tu(x,t) 
	:=
	\left \{ 
	\begin{aligned}
		& u(x,t) 
			&& \mbox{ for } (x,t)\in D,
\\		& \liminf_{x\uparrow \xi(t)} u(x,t) 
			&& \mbox{ for } (x,t)\in \Gamma,
	\end{aligned} 
	\right.
\]
then   the singularity at $x=\xi(t)$ is 
removable if and only if $\tu$ is continuously differentiable
 at $x=\xi(t)$ for any $t\in(0,T)$.


Next, we consider a higher dimensional singular set
whose spatial codimension is greater than or equal to $2$.
We reformulate our problem as follows. Let
$m\geq 1$, $N \geq m+2$, $T>0$ and 
$\bs=(s_1,s_2,\ldots,s_m) \in \R^m$.
We assume that the mapping
\[
	\xi(\bs,t) = (\xi^1(\bs,t),\xi^2(\bs,t),\ldots,\xi^{N}(\bs,t))
		: [0,1]^m \times [0,T] \rightarrow \R^N 
\]
is continuously differentiable with respect to $\bs$ 
and H\"older continuous with exponent $1/2$
with respect to $t$. 
Also, we assume that the Jacobian matrix of $\xi$ with
respect to $\bs$ is non-singular, that is, 
\begin{equation}\label{eq:rank}
\begin{aligned}
	\rank
		\begin{pmatrix}
			\xi_{s_1}^1(\sm,t)	&	\cdots	&	\xi_{s_m}^{1}(\sm,t)
		\\	\vdots				&	\ddots	&	\vdots
		\\	\xi_{s_1}^N(\sm,t)	&	\cdots	&	\xi_{s_m}^N(\sm,t)
		\end{pmatrix}
	= m
\end{aligned}
\end{equation}
for any $(\sm)\in [0,1]^m$ and $t\in[0,T]$.
We denote the singular set by
\[
	\Xi(t):= \{ \xi(\bs,t)\, : \, \bs \in [0,1]^m \}
\]
and define $\Gamma \subset \R^{N+1}$ by
\[
	\Gamma
	:= 
	\big \{ 
		(x,t) \in \R^{N+1} \, : \, x\in \Xi(t),
		\; t\in(0,T) 
	\big \}.
\]
We also define a distant between $x$ and $\Xi(t)$ by
\[
	d(x,\Xi(t)) 
	:=\min_{\bs\in[0,1]^m} |x-\xi(\bs,t)|.
\]
Furthermore, let $\Omega \subset \R^N$ be a domain such that
\[
 \Omega \supset \bigcup_{\ t\in[0,T]}\Xi (t) ,
\]
and define a domain $D \subset \R^{N+1}$ by 
\[
	D
	:= \big \{
		 (x,t) \in \R^{N+1} \, : \, 
		x\in \Omega \setminus 
			\Xi (t) 
			, \; t\in(0,T) 
		\big\}.
\]
Now we define removability of a higher dimensional singular set as follows.
For a solution of (\ref{eq:main}), 
 the singular set $\Xi(t)$ is said to be 
 removable if there exists a function $\tu$ which satisfies
the heat equation in $\Omega \times (0,T)$ 
in the classical sense and $u \equiv \tu$ on $D$.

Our results for higher dimensional singular sets are as follows.

\begin{thm}\label{th:removem} \ 
Let $N\geq m+3$. Suppose that $\xi$ satisfies $(\ref{eq:rank})$ 
and that $u$ satisfies $(\ref{eq:main})$ 
 in the classical sense.
Then the singular set $\Xi(t)$ is removable 
if and only if for any $0<t_1<t_2<T$ and $0<\eps <1$ 
there exists $0<r<1$ depending on $t_1,t_2,\eps$ such that
\begin{equation}\label{eq:epsm}
	|u(x,t)| 
	\leq \frac{\eps}{d(x,\Xi(t))^{N-m-2}},
		\quad 0<d(x,\Xi(t))<r
\end{equation}
for any $t \in [t_1, t_2]$.
\end{thm}

\begin{thm}\label{th:2removem} \ 
Let $N=m+2$. Suppose that $\xi$ satisfies $(\ref{eq:rank})$ 
and that $u$ satisfies $(\ref{eq:main})$ 
 in the classical sense.
Then the singular set $\Xi(t)$ is removable 
if and only if for any $0<t_1<t_2<T$ and $0<\eps <1$ 
the function $u$ satisfies
\[
	|u(x,t)| \leq 
	\eps \log \frac{1}{d(x,\Xi(t))},
		\quad 0<d(x,\Xi(t))<\eps
\]
for any $t \in [t_1, t_2]$.
\end{thm}

By an analogous method to Section \ref{sect:higher}, 
we can extend Theorem \ref{th:removem} to the case where
the singular set consists of 
$\Xi_1,\Xi_2,\ldots,\Xi_k$, each of which 
satisfies (\ref{eq:rank}) and may intersect with others. 
By regarding $\Xi_1,\Xi_2,\ldots,\Xi_k$ as local coordinates, 
the above theorems give a condition for the removability
in the case where 
the singular set is a compact $m$-dimensional $C^1$-manifold in $\R^N$.


Next, we show the existence of a solution of (\ref{eq:main})
whose singularity
moves in time and is not removable.
Again, let $N\geq 2$, $T>0$, and
$\Gamma\subset\R^{N+1}$ be defined as above.

The next result implies that the conditions in Theorem~\ref{th:remove0} 
and Theorem~\ref{th:2remove0} for
the removability are optimal in some sense.

\begin{thm}\label{th:example} \ 
Given any H\"older continuous function $\xi(t) : [0,T] \to \R^N$ 
with exponent $\alpha>1/2$, 
there exists $u$ defined on a 
neighborhood of $\Gamma$ such that $u$ satisfies 
$(\ref{eq:main})$ in the classical sense
but the singularity of $u$ at $x=\xi(t)$ is not removable.
\end{thm}

In Section \ref{sect:nonremove}, we give 
an example of a non-removable moving singularity.
In fact, this theorem will be proved by solving the following problem:
\begin{equation} \label{eq:atdelta}
	u_t - \Delta u = \delta(x - \xi(t)) 
		\quad \mbox{ in } \R^N \times (0,T), 
\end{equation}
where $\delta(\cdot)$ denote the Dirac distribution concentrated
at the point $0 \in \R^N$.
In this case, we can show that 
the singularity at $x=\xi(t)$ persists
for $t\in(0,T)$ and the solution satisfies
\[
\begin{aligned}
	& 
	u(x,t) 
	= \frac{1}{N(N-2)\omega_N} |x - \xi(t)|^{2-N} 
		+ o(|x - \xi(t)|^{2-N})
	&& \quad \mbox{if} \quad N\geq 3,
	\\ & 
	u(x,t) 
	= \frac{1}{2\pi} \log \Big( \frac{1}{|x-\xi(t)|} \Big) 
		+ o \Big( \log \frac{1}{|x-\xi(t)|} \Big)
	&& \quad \mbox{if} \quad N=2
\end{aligned}
\]
at $x=\xi(t)$, where we denote by $\omega_N$ 
volume of unit ball in $\R^N$.

This paper is organized as follows.
In Section \ref{sect:point} we prove 
Theorems \ref{th:remove0} and \ref{th:2remove0} 
by cutting a neighborhood of the singularity.
In Section \ref{sect:higher} we apply this method
to a higher dimensional singular set.
Section \ref{sect:nonremove} is devoted to 
the analysis of (\ref{eq:atdelta}).

\Sect{Removability of a moving singularity}\label{sect:point}
In this section, we consider removability of a moving singularity.
To show Theorem~\ref{th:remove0}, we give the following lemma.


\begin{lem}\label{lem:cut} \ 
Let $r>0$. Suppose that $\xi(t): [0,T]\rightarrow \R^N$ is 
H\"older continuous with exponent $\alpha>0$.
Then there exists a family of cut-off functions 
$\{\eta^r\}_{r>0}\subset C^\infty(\R^N \times (0,T))$ such that 
\[
	\eta^r(x,t)=
	\begin{cases}
	1 \quad & \mbox{ {\rm if} } |x-\xi(t)| > r,
	\\ 0 \quad & \mbox{ {\rm if} } |x-\xi(t)| < r/2,
	\end{cases}
\]
and
\[
	0\leq \eta^r \leq 1,\quad
	|\nabla \eta^r | \leq Cr^{-1},\quad
	|\Delta \eta^r | \leq Cr^{-2},\quad
	|(\eta^r)_t | \leq Cr^{-1/\alpha},
\]
where $C>0$ is a constant independent of $x,t$ and $r$.
\end{lem}


\begin{Proof}
Let $r>0$ be fixed. We take standard mollifier 
$\rho\in C^\infty(\R)$ by
\[
	\rho(t) :=
	\begin{cases}
	A e^{-1/(1-t^2)}
		\quad &\mbox{ if } |t|<1,
	\\ 0 
		\quad &\mbox{ if } |t|\geq1,
	\end{cases}
\]
where the constant $A>0$ is taken so that $\int_\R \rho(t)\,dt =1$. 
In addition, for each $\eps>0$, we set 
$\rho^{\eps}(t) := (1/\eps) \rho (t/\eps)$.
We express $\xi=(\xi_1, \xi_2, \ldots, \xi_N)$ and define
$\xi^\eps=(\xi_1^\eps, \xi_2^\eps, \ldots, \xi_N^\eps)$ by
\[
	\xi_i^\eps(t) :=
	\int_\R \rho^\eps (t-s)\xi_i(s) ds, \qquad i=1, 2, \ldots,N.
\]

Let $\eps>0$. By H\"older continuity of $\xi$, we obtain
\begin{equation}\label{eq:xieps}
	|\xi_i(t)-\xi_i^{\eps}(t)|
	\leq
	L \eps^\alpha
\end{equation}
for every $t\in[0,T]$, where $L>0$ is a H\"older constant.
Moreover, by changing variable 
$\tau=(t-s)/\eps$ and simple calculation, 
\[
\begin{aligned}
	(\xi^\eps _i)_t &=
	\frac{A}{\eps^2}
	\int_{t-\eps}^{t+\eps}
		\frac{-2(t-s)/\eps}{(1-((t-s)/\eps)^2)^2}
		\exp \Big( -\frac{1}{1- ((t-s)/\eps)^2} \Big)
		\xi_i(s)
	ds
	\\ & =
	\frac{A}{\eps}
	\int_{-1}^1
		\frac{-2\tau}{(1-\tau^2)^2}
		e^{-1/(1- \tau^2)}
		\xi_i(t-\eps \tau)
	d\tau.
\end{aligned}
\]
We remark that
\[
	\int_{-1}^1
		\frac{-2\tau}{(1-\tau^2)^2}
		e^{-1/(1- \tau^2)}
		\xi_i(t)
	d\tau=0.
\]
Then, by H\"older continuity, we have
\begin{equation}\label{eq:xiepst}
\begin{aligned}
	|(\xi^\eps _i)_t| &=
	\frac{A}{\eps}
	\Big|
	\int_{-1}^1
		\frac{-2\tau}{(1-\tau^2)^2}
		e^{-1/(1- \tau^2)}
		(
			\xi_i(t-\eps \tau) - \xi_i(t)
		)
	\,d\tau
	\Big|
	\\ & \leq
	AL \eps^{\alpha-1}
	\int_{-1}^1
		\frac{2|\tau|}{(1-\tau^2)^2}
		e^{-1/(1- \tau^2)}
	\,d\tau
	\\ & =
	2AL\eps^{\alpha-1}e^{-1} 
	\leq 
	AL\eps^{\alpha-1}.
\end{aligned}
\end{equation}

Now we define $\eta^r \in C^\infty(\R \times(0,T))$ by
\[
	\eta^r(x,t)
	:=
	\left\{
	\begin{aligned}
	&\frac{e^{-1/\sigma}}{e^{-1/\sigma} + e^{-1/(1-\sigma)}}
	&&\quad \mbox{ if } 
		\frac{7}{10}r < |x-\xi^\eps (t)| <\frac{4}{5}r,
	\\ &1
	&&\quad \mbox{ if }
		|x-\xi^\eps (t)| \geq \frac{4}{5}r,
	\\ &0
	&&\quad \mbox{ if }
		|x-\xi^\eps (t)| \leq \frac{7}{10}r,
	\end{aligned}
	\right.
\]
where
\[
	\sigma = \sigma(x,t;r) = 
	\frac{10}{r} 
	\Big(
		|x-\xi^\eps (t)| - \frac{7}{10}r
	\Big).
\]
It is clear that 
$0\leq \eta^r(x,t) \leq 1$.

Next, we take $\eps_r=(r/10NL)^{1/\alpha}$. 
By (\ref{eq:xieps}), we have
\[
\begin{aligned}
	|\xi(t)-\xi^{\eps_r}(t)| 
	& \leq
	|\xi_i(t)-\xi_i^{\eps_r}(t)|
		+ \cdots + |\xi_N(t)-\xi_N^{\eps_r}(t)|
	\\ & \leq
	L (\eps_r)^\alpha N = r/10.
\end{aligned}
\]
Here, $\eta^r(x,t)=1$  if   $|x-\xi(t)|> r$, because
\[
	|x-\xi^{\eps_r} (t)| \geq |x-\xi(t)| - |\xi(t)-\xi^{\eps_r} (t)|
	> r - (r/10) > 4r/5,
\]
and
 $\eta^r(x,t)=0$  if $|x-\xi(t)| < r/2$, because
\[
	|x-\xi^{\eps_r} (t)| \leq |x-\xi(t)| + |\xi(t)-\xi^{\eps_r} (t)|
	< (r/2) + (r/10) < 7r/10.
\]

Finally, we estimate first and second derivatives of $\eta^r$.
It suffices to calculate in the case where 
$7r/10 < |x-\xi^{\eps_r} (t)| <4r/5$.
In this case, we have $0<\sigma(x,t;r)<1$.
By direct calculation, we have
\[
	\nabla_x (\eta^r) = \frac{10}{r} X(\sigma) 
		\frac{x-\xi^{\eps_r} (t)}{|x-\xi^{\eps_r} (t)|},
	\qquad
	(\eta^r)_t = -\frac{10}{r} X(\sigma) 
		\frac
		{(x-\xi^{\eps_r} (t)) \cdot \xi^{\eps_r} _t(t)}
		{|x-\xi^{\eps_r} (t)|},
\]
where
\[
	X(\sigma)
	:=
	\frac{ e^{-1/\sigma} e^{-1/(1-\sigma)} }
		{ ( e^{-1/\sigma} + e^{-1/(1-\sigma)} )^2 }
	\Big(
		\frac{1}{\sigma^2} 
		+ \frac{1}{(1-\sigma)^2}
	\Big),
\]
and
\[
	\Delta_x (\eta^r)
	=
	\frac{100}{r^2} Y(\sigma),
\]
where
\[
\begin{aligned}
	Y(\sigma):=
	&
		\frac{ e^{-1/\sigma} e^{-1/(1-\sigma)} }
		{ ( e^{-1/\sigma} + e^{-1/(1-\sigma)} )^2 }
	\left[
		\Big(
			\frac{N-1}{\sigma+7}
		\Big)
		\Big(
			\frac{1}{\sigma^2}+\frac{1}{(1-\sigma)^2}
		\Big)
	\right.
	+
		\Big(
			\frac{1}{\sigma^4}+\frac{1}{(1-\sigma)^4}
		\Big)
		(1-2\sigma)
	\\ &\quad -
	\left.
		\frac{2}{ e^{-1/\sigma} + e^{-1/(1-\sigma)} }
		\Big(
			\frac{ e^{-1/\sigma} }{\sigma^4}
			+ \frac{ e^{-1/\sigma} - e^{-1/(1-\sigma)} }
				{\sigma^2 (1-\sigma^2)} 
			- \frac{ e^{-1/(1-\sigma)} }{(1-\sigma)^4}
		\Big)
	\right].
\end{aligned}
\]
Since $X(\sigma)$ and $Y(\sigma)$ belong to 
$C^\infty(0,1)$ and satisfy
\[
	\lim_{\sigma \downarrow 0} |X(\sigma)| 
	= \lim_{\sigma \uparrow 1} |X(\sigma)| 
	= \lim_{\sigma \downarrow 0} |Y(\sigma)| 
	= \lim_{\sigma \uparrow 1} |Y(\sigma)| = 0,
\]
we see that $X(\sigma)$ and $Y(\sigma)$ are bounded for $\sigma\in(0,1)$.
Moreover, by (\ref{eq:xiepst}), we obtain 
\[
	| (\eta^r)_t | 
	\leq 
	C_1 r^{-1} A L (\eps_r)^{\alpha-1} N
	=
	C_2 r^{-1/\alpha},
\]
where $C_1,C_2>0$ are constants independent of $x,t,r$.
Hence there exists a constant $C_3>0$ 
independent of $x,t,r$ such that
\[
	|\nabla \eta^r | \leq C_3 r^{-1}, \quad
	|\Delta \eta^r | \leq C_3 r^{-2}, \quad 
	|(\eta^r)_t | \leq C_3 r^{-1/\alpha}.
\]
The proof is complete.
\qed
\end{Proof}

\ 

\noindent
{\bf Proof of Theorem \ref{th:remove0}.} \ 
Necessity is easily proved by the same argument as in Section 3 of \cite{Hsu}.
Indeed, if the singularity of $u$ at $x=\xi(t)$ is removable, 
then $u$ is bounded near $x=\xi(t)$.

We prove sufficiency. 
Let $0<t_1<t_2<T$ and $0<\eps<1$.
By our assumption, there exists $r=r(t_1,t_2,\eps)>0$ 
such that (\ref{eq:eps0}) holds.
For each $t\in(0,T)$, we take any sequence 
$\{x_i(t)\}_{i=1}^\infty \subset \Omega\setminus\{\xi(t)\}$ 
such that $|x_i(t) -\xi(t)|\to 0$ as $i\to\infty$, and
set
\[
	\tu(x,t) 
	:=
	\left \{ 
	\begin{aligned}
		& u(x,t) 
			&& \mbox{ for } (x,t)\in D,
\\		& \liminf_{i\to \infty} u(x_i(t),t) 
			&& \mbox{ for } (x,t)\in \Gamma.
	\end{aligned} 
	\right.
\]
Our goal is to prove that $\tu$ satisfies the heat equation 
in $\Omega \times(0,T)$ in the classical sense.

First, we show 
$\tu \in L^1_{\loc}(\Omega \times(0,T))$.
For each $t\in[t_1, t_2]$, we denote
\[
	B(\xi(t),r) 
	:= \{ x\in \R^N \,:\, |x-\xi(t)|<r \}.
\]
By $N$-dimensional polar coordinates centered at $\xi(t)$, we have
\begin{equation}\label{eq:B}
	\int_{t_1}^{t_2}\int_{B(\xi(t),r)}
		|x-\xi(t)|^{2-N} \,dxdt
	= C_1(t_2-t_1)r^2
\end{equation}
for some $C_1=C_1(N)>0$.
Let $K\subset \R^N$ be a compact subset of $\Omega$.
Since $\xi(t) \in \Omega$ for $t \in [0,T]$, 
we can take $r=r(t_1,t_2,\eps)>0$ so small that 
 $B(\xi(t),r) \subset \Omega$ for every $t\in[t_1, t_2]$.
By (\ref{eq:eps0}) and (\ref{eq:B}), there exists $C_2>0$ such that
\[
\begin{aligned}
	\int_{t_1}^{t_2}\int_K
		|\tu(x,t)| \,dxdt 
	\leq & 
	\int_{t_1}^{t_2} \int_{K\setminus B(\xi(t,r))}
		 |u(x,t)| \,dxdt
	+
	\eps \int_{t_1}^{t_2}\int_{B(\xi(t),r)}
		|x-\xi(t)|^{2-N} \,dxdt
	\\ \leq & 
		C_2 +\eps C_1 (t_2-t_1)r^2 
		<\infty.
	\end{aligned}
\]
Since $0<t_1<t_2<T$ are arbitrary, 
we have $\tu\in L^1_{\loc}(\Omega \times (0,T))$.

Next, we show that 
$\tu$ satisfies the heat equation
in $\Omega \times(0,T)$ in the distribution sense.
For this purpose, 
we need a family of cut-off functions 
$\{\eta^r\}_{r>0} \subset C^{\infty}(\R^N \times (0,T))$ 
such that
\[
	\eta^r(x,t) =
	\left \{
		\begin{aligned}
		& 0 && \mbox{ if }\quad |x-\xi(t)| < r/2,
\\		& 1 && \mbox{ if }\quad |x-\xi(t)| > r,
		\end{aligned}
	\right.
\]
and
\begin{equation}\label{eq:eta}
	0\leq \eta^r \leq 1,\quad
	|\nabla \eta^r| \leq C_3 r^{-1},\quad
	|\Delta \eta^r| \leq C_3 r^{-2},\quad
	|(\eta^r)_t| \leq C_3 r^{-2},
\end{equation}
where $C_3>0$ is a constant independent of $x,t$ and $r$.
By Lemma \ref{lem:cut} and the assumption that 
$\xi(t)$ is H\"older continuous with exponent $1/2$,
we can take such $\{\eta^r\}$.
Now, let $\phi \in C^{\infty}_0(\Omega \times (0,T))$ 
be a test function. 
Since $\eta^r \phi \in C^{\infty}_{0}(\Omega \times (0,T))$, and 
$\tu$ is a classical solution of (\ref{eq:main}), we have
\begin{equation}\label{eq:kdist}
	\begin{aligned}
		 \int_{\Omega}&\{ \tu(x,t_2) \phi(x,t_2)\eta^r(x,t_2)
			-\tu(x,t_1)\phi(x,t_1)\eta^r(x,t_1) \}\,dx
\\		& =
		 \int_{t_1}^{t_2} \int_{\Omega}
		 	\tu \{ (\phi \eta^r)_t
		+ \Delta(\phi\eta^r) \} \,dx dt.
	\end{aligned}
\end{equation}

Here, we claim that the following convergence properties hold:
\begin{align}
& \limsup_{\eps\to0}
	\left|
		\int_{t_1}^{t_2}\int_{\Omega}\tu \Delta\phi \,dxdt
		- \int_{t_1}^{t_2}\int_{\Omega}\tu\Delta(\phi \eta^r) \,dxdt
	\right|
	= 0, \label{eq:conv1}
\\ & \limsup_{\eps\to0}
	\left|
		\int_{t_1}^{t_2}\int_{\Omega}\tu \phi_t \,dxdt
		- \int_{t_1}^{t_2}\int_{\Omega}\tu(\phi \eta^r)_t \,dxdt
	\right|
	= 0, \label{eq:conv2}
\\ & \limsup_{\eps\to0}
	\left|
		\int_{\Omega}\tu(x,t_1) \phi(x,t_1) \,dx
		- \int_{\Omega}\tu(x,t_1)\phi(x,t_1)\eta^r(x,t_1) \,dx
	\right|
	= 0, \label{eq:conv3}
\\ & \limsup_{\eps\to0}
	\left|
		\int_{\Omega}\tu(x,t_2) \phi(x,t_2) \,dx
		- \int_{\Omega}\tu(x,t_2)\phi(x,t_2)\eta^r(x,t_2) \,dx
	\right|
	= 0. \label{eq:conv4}
\end{align}
To show (\ref{eq:conv1}), we rewrite
\begin{equation}\label{eq:Ir}
	\begin{aligned}
		 \int_{t_1}^{t_2} & \int_{\Omega}
			\tu \Delta\phi \,dxdt
		- \int_{t_1}^{t_2}\int_{\Omega}
			\tu\Delta(\phi \eta^r) \,dxdt 
\\		&= \int_{t_1}^{t_2}\int_{\Omega}
			\tu(1-\eta^r)\Delta \phi \,dx dt
		-2\int_{t_1}^{t_2}\int_{\Omega}
			\tu \nabla \phi\cdot \nabla \eta^r \,dx dt
		-\int_{t_1}^{t_2}\int_{\Omega}
			\tu \phi\Delta \eta^r \,dx dt
\\		&=: \, I_{1,r} -2 I_{2,r} - I_{3,r}.
	\end{aligned}
\end{equation}
By (\ref{eq:eps0}) and (\ref{eq:eta}), 
for sufficiently small $r=r(t_1,t_2,\eps)>0$, 
we have the inequalities
\[
	\begin{aligned}
		|I_{1,r}| 
		& \leq 
		\| \Delta \phi \|_{L^{\infty}(\Omega\times(0,T))} \eps
		\int_{t_1}^{t_2}\int_{B(\xi(t),r)}
			|x-\xi(t)|^{2-N} \,dxdt,
\\		|I_{2,r}| 
		& \leq 
		\| \nabla \phi \|_{L^{\infty}(\Omega\times(0,T))}
		C_3 \frac{\eps}{r}
		\int_{t_1}^{t_2}\int_{B(\xi(t),r)}
			|x-\xi(t)|^{2-N} \,dxdt,
\\		|I_{3,r}| 
		& \leq 
		\| \phi \|_{L^{\infty}(\Omega\times(0,T))}
		C_3 \frac{\eps}{r^2}
		\int_{t_1}^{t_2}\int_{B(\xi(t),r)}
			|x-\xi(t)|^{2-N} \,dxdt.
	\end{aligned}
\]
Hence, by (\ref{eq:B}) and $r\in(0,1)$, we have
\[
	\begin{aligned}
		|I_{1,r}| 
		& \leq 
		\| \Delta \phi \|_{L^{\infty}(\Omega\times(0,T))} 
		C_1(t_2-t_1) \eps r^2
		 \leq 
		C_4 \eps,
\\		|I_{2,r}| 
		& \leq 
		\| \nabla \phi \|_{L^{\infty}(\Omega\times(0,T))} 
		C_1 C_3(t_2-t_1) \eps r
		 \leq
		C_4 \eps,
\\		|I_{3,r}| 
		& \leq 
		\| \phi \|_{L^{\infty}(\Omega\times(0,T))}
		C_1 C_3(t_2-t_1) \eps
		 \leq C_4 \eps
\end{aligned}
\]
for some $C_4>0$.
Hence we obtain (\ref{eq:conv1}).
Similarly we obtain (\ref{eq:conv2}), (\ref{eq:conv3}) and (\ref{eq:conv4})
from above estimates.

Thus, the function $\tu$ satisfies
\begin{equation}\label{eq:dist}
	\int_{\Omega}\{\tu(x,t_2)\phi(x,t_2)-\tu(x,t_1)\phi(x,t_1)\} \,dx
	=\int_{t_1}^{t_2}\int_{\Omega}\tu(\phi_t +\Delta \phi) \,dxdt
\end{equation}
for any $\phi \in C^{\infty}_{0}(\Omega \times (0,T))$.
Since $0<t_1<t_2<T$ be arbitrary, 
the function $\tu\in L^1_{\loc}(\Omega\times(0,T))$ satisfies 
the heat equation in $\Omega \times (0,T)$ in the distribution sense.
By using the Weyl lemma for the heat equation 
(see, e.g., Section 6 of \cite{G} or \cite{S}), 
$\tu$ satisfies the heat equation in $\Omega \times(0,T)$ 
in the classical sense.
Since $\tu=u$ in $D$, 
 the singularity of $u$ at $x=\xi(t)$ is removable. 
\Qed

\noindent
{\bf Proof of Theorem \ref{th:2remove0}.} \ 
We prove only sufficiency. 
Let $0<t_1<t_2<T$ and $0<\eps<1$, and define 
 $\tu$ as in the proof of Theorem \ref{th:remove0}.
By $2$-dimensional polar coordinates, we have
\begin{equation}\label{eq:2B}
	\begin{aligned}
		\int_{t_1}^{t_2}\int_{B(\xi(t),\eps)}
		\log \frac{1}{|x-\xi(t)|} \,dxdt
		&\leq C_1(t_2-t_1)
		\left(
			1+\log (1/\eps)
		\right) \eps^2
	\end{aligned}
\end{equation}
for some $C_1>0$.
This implies $\tu\in L^1_{\loc}(\Omega \times (0,T))$.

We show that $\tu$ satisfies the heat equation in $\Omega \times(0,T)$ 
in the distribution sense.
Let $\phi \in C_0^\infty(\Omega \times (0,T))$.
By Lemma \ref{lem:cut} and  the assumption that 
$\xi(t)$ is H\"older continuous with exponent $1/2$,
we can take $\{\eta^\eps\}_{\eps>0}\subset 
C^{\infty}(\R^N \times(0,T))$ such that
\[
	\eta^\eps(x,t)=
	\left\{ 
		\begin{aligned}
			& 0 && \mbox{ if }\quad |x-\xi(t)| < \eps/2,
		\\	& 1 && \mbox{ if }\quad |x-\xi(t)| > \eps,
		\end{aligned}
	\right.
\]
and
\begin{equation}\label{eq:2eta}
	0\leq \eta^\eps \leq 1,\quad
	|\nabla \eta^\eps| \leq C_2 \eps^{-1},\quad
	|\Delta \eta^\eps| \leq C_2 \eps^{-2},\quad
	|(\eta^\eps)_t| \leq C_2 \eps^{-2}
\end{equation}
for some $C_2>0$. 
Since $\tu$ satisfies (\ref{eq:main}), the equality
 (\ref{eq:kdist}) holds for $r=\eps$.
Again, we claim that the convergence properties 
(\ref{eq:conv1}), (\ref{eq:conv2}), (\ref{eq:conv3}) and 
(\ref{eq:conv4}) hold for $r=\eps$.
Let $I_{1,\eps},I_{2,\eps}$ 
and $I_{3,\eps}$ be defined as in (\ref{eq:Ir}) with $r=\eps$.
By (\ref{eq:2eps0}) and (\ref{eq:2eta}), 
for sufficiently small $\eps>0$, 
we have
\[
	\begin{aligned}
		|I_{1,\eps}| 
		& \leq 
		\| \Delta \phi \|_{L^{\infty}(\Omega\times(0,T))}
		\eps \int_{t_1}^{t_2}\int_{B(\xi(t),\eps)}
		\log \frac{1}{|x-\xi(t)|} \,dxdt,
		\\ |I_{2,\eps}| 
		& \leq 
		\| \nabla \phi \|_{L^{\infty}(\Omega\times(0,T))}
		C_2 \int_{t_1}^{t_2}\int_{B(\xi(t),\eps)}
		\log \frac{1}{|x-\xi(t)|} \,dxdt,
		\\ |I_{3,\eps}| 
		& \leq 
		\| \phi \|_{L^{\infty}(\Omega\times(0,T))}
		C_2 \frac{1}{\eps} 
		\int_{t_1}^{t_2}\int_{B(\xi(t),\eps)}
		\log \frac{1}{|x-\xi(t)|} \,dxdt.
\end{aligned}
\]
Hence by (\ref{eq:2B}), we have
\[
	\begin{aligned}
	|I_{1,\eps}|
		& \leq 
		\| \Delta \phi \|_{L^{\infty}(\Omega\times(0,T))} 
		C_1(t_2-t_1) \left( 1+\log (1/\eps) \right) \eps^3
		 \leq C_3 \eps \log (1/\eps),
	\\ 
	|I_{2,\eps}| 
		& \leq 
		\| \nabla \phi \|_{L^{\infty}(\Omega\times(0,T))} 
		C_1C_2(t_2-t_1)
		\left(
			1+\log (1/\eps)
		\right) \eps^2
		 \leq C_3 \eps \log (1/\eps),
	\\ 
	|I_{3,\eps}| 
		& \leq 
		\| \phi \|_{L^{\infty}(\Omega\times(0,T))}
		C_1C_2(t_2-t_1) 
		\left( 
			1+\log (1/\eps) 
		\right) \eps
		 \leq C_3 \eps \log (1/\eps)
\end{aligned}
\]
for some $C_3>0$.
Hence we obtain (\ref{eq:conv1}).
Similarly we obtain (\ref{eq:conv2}), (\ref{eq:conv3}) and (\ref{eq:conv4})
from above estimates.
These imply that $\tu\in L^1_{\loc}(\Omega\times(0,T))$ 
satisfies the heat equation
in $\Omega \times (0,T)$ in the distribution sense.
The remainder is the same as in the proof of Theorem \ref{th:remove0}.
\Qed


\Sect{Removability of a  singular set}\label{sect:higher}
Let 
 $\Xi(t)\subset\R^N$, $D\subset\R^{N+1}$, $\Gamma\subset\R^{N+1}$, 
and $\Omega\subset\R^N$ are the sets defined in Section \ref{sect:intro}.
To show Theorems \ref{th:removem} and \ref{th:2removem}, 
we give the following estimates.

\begin{lem}\label{lem:xi} \ 
There exists $C_1=C_1(N,m)>0$ and $C_2=C_2(m)>0$ such that
for every sufficiently small $r>0$, 
\begin{align}
	& \int_{A_{r,t}} d(x,\Xi(t))^{m+2-N} \,dx 
		\leq C_1 r^2
		\quad
	&& \mbox{if}\quad N\geq{m+3}, \label{eq:xi} 
	\\ & \int_{A_{r,t}} \log \frac{1}{d(x,\Xi(t))} \,dx 
		\leq C_2 r^2\left( 1+\log\frac{1}{r} \right)
		\quad
	&& \mbox{if}\quad N=m+2 \label{eq:2xi}
\end{align}
for any $t\in(0,T)$, 
where $A_{r,t} := \{ x\in \R^N : d(x,\Xi(t))<r\}$.
\end{lem}

\begin{Proof}
We prove the lemma only in the case $N\geq m+3$. 
In fact, (\ref{eq:2xi}) can be proved 
in the same manner as (\ref{eq:xi}).
Let $t\in(0,T)$ be fixed.
We extend the domain of the function $\xi$ to $[a,b]^m\times[0,T]$ 
with $a<0$ and $b>1$.
That is, we take a mapping
\[
	\txi(\bs,t) = (\txi^1(\bs,t),\txi^2(\bs,t),\ldots,\txi^{N}(\bs,t))
	:[a,b]^m \times [0,T] \rightarrow \R^N 
\] 
such that $\txi$ is continuously differentiable 
in $\bs$ and continuous in $t$.
In addition, we assume that $\txi$ satisfies (\ref{eq:rank}) and
\[
	\txi \big|_{[0,1]^m\times[0,T]} = \xi,
	\quad
	\txi^j_{s_i} \big|_{[0,1]^m\times[0,T]} = \xi^j_{s_i},
		\quad i=1,2,\ldots,m, \quad j=1,2,\ldots,N.
\]
We define
\[
	\tXi(t)
	:= \{ \txi(\bs,t)\, : \, \bs \in [a,b]^m \}.
\]

For each $\bs\in (a,b)^m$, 
let $\Pi_{r,t}(\bs)$ be a subset 
of a normal plane of $\tXi(t)$ at $\txi(\bs,t)$ given by 
\[
	\Pi_{r,t}(\bs)
	:= \{
		 x\in A_{r,t} \, : \,
			(x-\txi(\bs,t))\cdot \txi_{s_i}(\bs,t)=0
		 \;	\mbox{ for any } i=1,2,\ldots,m
		\}.
\]
Since $\txi(\cdot,t)$ is defined on a compact set, 
there exists a sufficiently small $r>0$ such that
\begin{equation}\label{eq:shortest}
	d(x,\tXi(t)) 
		= |x-\txi(t)|,\quad x\in \Pi_{r,t}(\bs)
\end{equation}
for each $\bs\in(a,b)^m$.
Again by compactness, we have
\begin{equation}\label{eq:m}
	M:=\max_{t\in[0,T]}
		\int_{ \tXi(t) } \,d\sigma^m
	<\infty,
\end{equation}
where $d\sigma^m$ is an $m$-dimensional surface element.
%
%
%
Since $\txi$ satisfies (\ref{eq:rank}), 
$\Pi_{r,t}(\bs)$ is an $(N-m)$-dimensional  subspace of $\R^N $.
Therefore, for
 each $\bs\in (a,b)^m$, there exists a congruent transformation
 $P_\bs:\R^N\rightarrow\R^N$  such that
\[
	P_\bs x
	= (y_1,y_2,\ldots,y_{N-m},0,\ldots,0),
	\quad
	x\in \Pi_{r,t}(\bs)
\]
for some $y_1,y_2,\ldots,y_{N-m}\in \R$.
Now, by using $(N-m)$-dimensional polar coordinates, we obtain
\begin{equation}\label{eq:N-m}
	\int_{P_\bs(\Pi_{r,t}(\bs))}
	|y|^{m+2-N} \,dy_1 dy_2 \cdots dy_{N-m}
	= C_3 \int_0^r \rho^{m+2-N+(N-m-1)} \,d\rho
	=C_4 r^2,
\end{equation}
where $C_3,\,C_4 >0$ depend on $N,m$ but not on $\bs, t$.

Recall that 
 the congruent transformations preserve a distance
between any two points and that the
function $\txi$ is an extension of $\xi$. 
Hence by choosing sufficiently small $r>0$ again if necessary, 
we have the estimate
\[
	\int_{A_{r,t}} d(x,\xi(\cdot,t))^{m+2-N} \,dx
	\leq MC_4 r^2
\]
by using (\ref{eq:shortest}), (\ref{eq:m}) and (\ref{eq:N-m}).
Thus we obtain (\ref{eq:xi}).
\Qed
\end{Proof}

\noindent
{\bf Proof of Theorem \ref{th:removem}.} \ 
We adopt the same approach as in the proof of Theorem \ref{th:remove0}, 
so we state the outline only.

Let $0<t_1<t_2<T$ and $0<\eps<1$.
By our assumption, there exists $r=r(t_1,t_2,\eps)>0$ 
such that (\ref{eq:epsm}) holds.
For $t\in(0,T)$, we take any sequence 
 $\{x_i(t)\}_{i=1}^{\infty} 
\subset \Omega \setminus \Xi(t)$ 
 such that $d(x_i(t), \xi(\cdot,t))\to 0$ as $i\to\infty$, 
and set
\[
	\tu(x,t) :=
	\left\{
		\begin{aligned}
			& u(x,t) 
			&& \mbox{ for } (x,t)\in D,
		\\	& \liminf_{i\to\infty} u(x_i(t),t) 
			&& \mbox{ for } (x,t)\in \Gamma.
		\end{aligned}
	\right.
\]
By Lemma \ref{lem:xi}, 
we obtain $\tu\in L^1_{\loc}(\Omega \times (0,T))$.
We show that $\tu$ satisfies (\ref{eq:main}) in $\Omega \times(0,T)$
in the distribution sense.
Let $\phi \in C_0^\infty(\Omega \times (0,T))$. 
By an argument similar to Lemma \ref{lem:cut}, 
we can take 
$\{\eta^r\}_{r>0}\subset C^\infty(\Omega \times(0,T))$ 
such that
\[
	\eta^r(x,t)=
	\left\{
		\begin{aligned}
			& 0 && \mbox{ if } 
				\quad d(x,\xi(\cdot,t)) < r/2,
		\\	& 1 && \mbox{ if } 
				\quad d(x,\xi(\cdot,t)) > r,
		\end{aligned}
	\right.
\]
and $\eta^r$ satisfies the condition (\ref{eq:eta}) 
for some $C>0$.
Since $\tu$ satisfies (\ref{eq:main}), we have (\ref{eq:kdist}).
By Lemma \ref{lem:xi} and an argument similar to
Section \ref{sect:point}, we obtain (\ref{eq:dist}).
That is, the function $\tu\in L^1_{\loc}(\Omega\times(0,T))$ 
satisfies the heat equation in $\Omega\times(0,T)$ 
 in the distribution sense.
The remainder is the same as in the proof of Theorem \ref{th:remove0}.
\Qed

Since (\ref{eq:2xi}) holds, we can show Theorem \ref{th:2removem} 
in the same way. We omit details of the proof.



\Sect{Non-removable singularity}\label{sect:nonremove}
In this section, we consider the case where 
a singularity move in time and is not removable.
Without loss of generality, we take $\Omega=\R^N$.
Let $N\geq 2$ and $T>0$. 
We assume that $\xi : [0,T]\rightarrow\R^N$ is 
arbitrarily given continuous function. 

To show Theorem \ref{th:example}, we 
solve the equation (\ref{eq:atdelta}). 
In this paper, we say that $u$ satisfies 
(\ref{eq:atdelta}) in the distribution sense 
if $u$ belongs to $L^1_{\loc}(\R^N\times(0,T))$ and satisfies
\begin{equation}\label{eq:deltadefi}
	\int_0^T\int_{\R^N}(-\phi_t-\Delta \phi)u \,dxdt
	= \int_0^T \phi(\xi(t),t) \,dt
\end{equation}
for any $\phi\in C_0^\infty(\R^N\times(0,T))$.
Now, we denote by 
\[
\Phi(x,t):=(4\pi t)^{-N/2}\exp (-|x|^2/4t)
\]
the fundamental solution of the heat equation.
Moreover, we define $F$ in $\R^N\times(0,T)$ by
\[
	F(x,t) :=
	\int_0^t \Phi(x-\xi(s),t-s) \,ds.
\]
In the following, we will show that $F$ satisfies 
(\ref{eq:atdelta}) in the distribution sense.
In addition, we will give upper and lower estimates of $F$, 
and we will see that $F$ is an example of 
Theorem \ref{th:example}.

\begin{pro}\label{pro:Fheat} \ 
The function $F$ satisfies $(\ref{eq:main})$ 
in the classical sense.
\end{pro}

To show Proposition \ref{pro:Fheat}, we give the following lemma.

\begin{lem}\label{lem:Fdelta} \ 
The function $F$ satisfies $(\ref{eq:atdelta})$ 
in the distribution sense.
\end{lem}

\begin{Proof}
First, we show $F\in L_{\loc}^1(\R^N\times(0,T))$.
By simple calculation, we have
\[
	\begin{aligned}
		\int_0^T \int_{\R^N} F(x,t) \,dxdt
		&= \int_0^T \int_0^t 
		\left(
			\int_{\R^N}\Phi(x-\xi(s),t-s) \,dx
		\right) \,dsdt
		\\&= \int_0^T \int_0^t \,dsdt 
		= \frac{1}{2}T^2 < \infty,
	\end{aligned}
\]
so that $F\in L^1(\R^N\times(0,T))$. 
In particular, $F$ belongs to $L^1_{\loc}(\R^N\times(0,T))$. 

Next, we show that $F$ satisfies (\ref{eq:deltadefi}).
For this purpose, 
let $\phi\in C_0^\infty(\R^N\times(0,T))$ be a test function.
For each $t\in(0,\tau)$, we take $\tau\in(0,t)$ and 
define $F^\tau$ by 
\[
	F^\tau(x,t)
	=\int_0^{t-\tau} \Phi(x-\xi(s),t-s) \,ds.
\]
Here $F^\tau$ is bounded for each fixed $\tau$, that is, 
there exist $C_1(N),\,C_2(N)>0$ such that
\[
	0
	\leq F^\tau(x,t) 
	\leq C_1(N)\int_0^{t-\tau}(t-s)^{-N/2} \,ds
	\leq C_2(N) \tau^{(2-N)/2}
\]
for each $t\in(0,T)$.
Then, integrating by parts yields
\[
	\begin{aligned}
		\int_0^T & \int_{\R^N}
			(-\phi_t -\Delta\phi) F^\tau \,dxdt
		\\ &=
		\int_0^T\int_{\R^N}
			(-\phi_t -\Delta\phi)
				\left(
					\int_0^{t-\tau} \Phi(x-\xi(s),t-s) \,ds
				\right)
		\,dxdt
		\\ &= 
		\int_0^T\int_{\R^N}
			\phi(x,t) \Phi(x-\xi(t-\tau),\tau) \,dxdt
		\\ &\quad + 
		\int_0^T\int_{\R^N}\phi(x,t)
			\left( 
				\int_0^{t-\tau} 
				\{
					\Phi_t(x-\xi(s),t-s)
					-\Delta \Phi(x-\xi(s),t-s) 
				\} \,ds
			\right) \,dxdt
		\\ & =
		\int_0^T\int_{\R^N}
			\phi(x,t)\Phi(x-\xi(t-\tau),\tau) \,dxdt.
	\end{aligned}
\]
Similarly from Section 2.3.1 of \cite{E}, we see that
\begin{equation}\label{eq:right}
	\lim_{\tau \to 0} \int_{\R^N}
		\phi(x,t)\Phi(x-\xi(t-\tau),\tau) \,dx 
	=
	\phi(\xi(t),t)
\end{equation}
for each $t\in(0,T)$.

For the reader's convenience, we give a proof of (\ref{eq:right}). 
Let $0<t<T$ and $\eps>0$ be fixed.
We choose $\delta>0$ such that
\begin{equation}\label{eq:eps}
	|\phi(x,t)-\phi(\xi(t),t)|<\eps
\end{equation}
for any $|x-\xi(t)|<\delta$.
Then, we have
\[
	\begin{aligned}
		\bigg|
			\int_{\R^N} & 
				 \phi(x,t) \Phi(x-\xi(t-\tau),\tau) \,dx
			-\phi(\xi(t),t)
		\bigg|
		\\ \leq & 
			\int_{\R^N} |\phi(x,t)-\phi(\xi(t),t)|
				\Phi(x-\xi(t-\tau),\tau) \,dx
		\\ = & 
			\int_{B(\xi(t),\delta)}
				+ \int_{\R^N\setminus B(\xi(t),\delta)} 
		=: I_1 + I_2.
	\end{aligned}
\]
First, 
by $(\ref{eq:eps})$, we have an estimate of $I_1$ as
\[
	I_1 
	\leq \eps \int_{\R^N} \Phi(x-\xi(t-\tau),\tau) \,dx 
	= \eps.
\]
Next, we give an estimate of $I_2$.
If $|x-\xi(t)|\geq\delta$ and 
$|\xi(t)-\xi(t-\tau)|\leq \delta/2$, then
\[
	|x-\xi(t)|
	\leq |x-\xi(t-\tau)| + |\xi(t-\tau)-\xi(t)|
	\leq |x-\xi(t-\tau)| + \frac{1}{2}|x-\xi(t)|
\]
Hence $|x-\xi(t-\tau)|\geq |x-\xi(t)|/2$.
By simple calculation, 
\[
\begin{aligned}
	I_2 & \leq 
	2 \| \phi \|_{L^\infty(\R^N\times(0,T))}
	\int_{\R^N \setminus B(\xi(t),\delta)}
		(4\pi\tau)^{-N/2}
		\exp
			\left(
				-\frac{|x-\xi(t-\tau)|^2}{4\tau} 
			\right) \,dx
	\\ & \leq 
	C_3\tau^{-N/2}\int_{\R^N \setminus B(\xi(t),\delta)}
		\exp
			\left(
				-\frac{|x-\xi(t)|^2}{16\tau} 
			\right) \,dx
	\\ &= 
	C_4\tau^{-N/2}\int_\delta^\infty
	r^{N-1}
		\exp 
			\left(
				-\frac{r^2}{16\tau} 
			\right) \,dr
	\\ &= 
	C_5\int_{\delta / 4\sqrt{\tau}}^\infty
	\sigma^{N-1}e^{-\sigma^2} \,d\sigma
	\to 0 \quad \mbox{as} \quad \tau\to0,
\end{aligned}
\]
where $C_3,\,C_4,\,C_5>0$ are constants independent of $\tau$, 
and $r=4\sqrt{\tau}\sigma$.
Therefore, if we have $|\xi(t)-\xi(t-\tau)|\leq \delta/2$ 
and take $\tau>0$ is sufficiently small, 
then we obtain $I_1+I_2\leq \eps$.
Thus it is shown that (\ref{eq:right}) holds.

From (\ref{eq:right}) and the Lebesgue theorem, 
we see that $F$ satisfies (\ref{eq:deltadefi}), 
that is, 
\begin{equation}\label{eq:Fsol}
	\int_0^T\int_{\R^N}
		(-\phi_t-\Delta \phi) F(x,t) \,dxdt 
	= \int_0^T \phi(\xi(t),t) \,dt.
\end{equation}
Hence the function $F$ satisfies 
(\ref{eq:atdelta}) in the distribution sense.
\Qed
\end{Proof}

\noindent
{\bf Proof of Proposition \ref{pro:Fheat}.} \ 
Let $\psi \in C_0^\infty(D)$ be a test function, in
 particular, $\psi\in C_0^\infty(\R^N\times(0,T))$.
By (\ref{eq:Fsol}), we have
\[
	\int_0^T\int_{\R^N}
		(-\psi_t-\Delta \psi)F \,dxdt
	= \int_0^T \psi(\xi(t),t) \,dt.
\]
Since $\psi(\xi(t),t)=0$ for any $t\in(0,T)$, 
we obtain
\[
	\int_0^T\int_{\R^N}
		(-\psi_t-\Delta \psi)F \,dxdt 
	=0.
\]
Hence $F\in L^1{(\R^N\times(0,T))}$ satisfies the heat equation 
in $D$ in the distribution sense.
By the Weyl lemma for the heat equation, 
we conclude that $F$ satisfies (\ref{eq:main}) 
in the classical sense.
\Qed

\begin{pro}\label{pro:uplow} \ 
Let $N\geq3$.
Suppose that $\xi$ is H\"older continuous
with exponent $\alpha >1/2$.
Then for each $t\in(0,T)$ the function $F(x,t)$ satisfies
\[
	F(x,t) = 
	\frac{1}{N(N-2)\omega_N} |x-\xi(t)|^{2-N}
		+ o(|x-\xi(t)|^{2-N}) 
	\quad \mbox{ as } x\to\xi(t),
\]
where $\omega_N$ is the volume of unit ball in $\R^N$.
\end{pro}

\begin{Proof}
We fix  $t\in(0,T)$ and  set $z:=x-\xi(t)$. 
By changing variable $t-s=|z|^2/(4\sigma)$, 
we have
\begin{equation}\label{eq:change}
	\begin{aligned}
		F(x,t) 
		&= 
		\int_0^t (4\pi(t-s))^{-N/2}
			\exp
			\left(
				-\frac{|z+\xi(t)-\xi(s)|^2}{4(t-s)} 
			\right) \,ds
		\\&=
		4^{-1}\pi^{-N/2} |z|^{2-N} 
		\int_{|z|^2/4t}^\infty
			\sigma^{(N/2)-2}
				\exp
				\left( - \left| 
					\sigma^{1/2} \frac{z}{|z|}
					+ \frac{1}{2}\frac{\xi(t)-\xi(s)}{(t-s)^{1/2}}
				\right|^2 \right) \,d\sigma
		\\&=:
		4^{-1}\pi^{-N/2} |z|^{2-N} I(z,t).
	\end{aligned}
\end{equation}
Here, we rewrite $I(z,t)$ as
\[
\begin{aligned}
	I(z,t) =
	\int_0^\infty 
			\sigma^{(N/2)-2} &  e^{-\sigma}
				\exp
				\left( 
					- \sigma^{1/2} \frac{z}{|z|}
					\cdot \frac{\xi(t)-\xi(s)}{(t-s)^{1/2}} 
				\right)
	\\ & 
				\exp
				\left( 
					- \frac{1}{4} \frac{|\xi(t)-\xi(s)|^2}{(t-s)} 
				\right)
				\chi_{[\, |z|^2/4t, \infty)}(\sigma)
	\, d\sigma,
\end{aligned}
\]
where $\chi_A$ is a indicator function of $A$.

In order to apply the Lebesgue theorem to $I(z,t)$, 
we construct a dominating integrable function as follows.
By H\"older continuity of $\xi$, 
for sufficiently small $|z|>0$, we have
\[
\begin{aligned}
	\sigma^{(N/2)-2} &  e^{-\sigma}
		\exp
		\left( 
			- \sigma^{1/2} \frac{z}{|z|}
			\cdot \frac{\xi(t)-\xi(s)}{(t-s)^{1/2}} 
		\right)
		\exp
		\left( 
			- \frac{1}{4} \frac{|\xi(t)-\xi(s)|^2}{(t-s)} 
		\right)
		\chi_{[\, |z|^2/4t, \infty)}(\sigma)
	\\ & \leq
	\sigma^{(N/2)-2} e^{-\sigma} 
		\exp
		\left(
			L \sigma^{1/2}
			\Big( \frac{|z|^2}{4\sigma} \Big)^{\alpha-(1/2)}
		\right)
	\leq
	\sigma^{(N/2)-2} e^{-\sigma+\sigma^{1-\alpha}},
\end{aligned}
\]
where $L>0$ is a H\"older constant.
Since $\alpha > 1/2$,  we see that
$\sigma^{(N/2)-2} e^{-\sigma+\sigma^{1-\alpha}}$ 
becomes a dominating integrable function.
On the other hand, by using H\"older continuity of $\xi$ again, 
we have
\[
\begin{aligned}
	&\Big| 
		- \sigma^{1/2} \frac{z}{|z|}
		\cdot \frac{\xi(t)-\xi(s)}{(t-s)^{1/2}} 
	\Big|
	\leq
	\frac{L}{4^{\alpha-(1/2)}} \sigma^{1-\alpha} 
		|z|^{2\alpha-1}
	\to 0
	\quad \mbox{ as } |z|\to0,
	\\ &
	\Big|
		- \frac{1}{4} \frac{|\xi(t)-\xi(s)|^2}{(t-s)} 
	\Big|
	\leq
	\frac{L^2}{4^{2\alpha}} \sigma^{-2\alpha+1}
		|z|^{4\alpha-2}
	\to 0
	\quad \mbox{ as } |z|\to0
\end{aligned}
\]
for each $\sigma\in(0,\infty)$. 
Hence by the Lebesgue theorem,  we obtain
\[
	\lim_{|z|\to 0} I(z,t) 
	=
	 \int_0^\infty 
	 	\sigma^{(N/2)-2} e^{-\sigma}
	\,d\sigma
	=
	\Gamma \Big( \frac{N}{2} - 1 \Big) 
	=
	\frac{ 4 \pi^{N/2}}{ N (N-2) \omega_N},
\]
where $\Gamma$ denotes the gamma function.
Hence by (\ref{eq:change}),
we obtain
\[
	\lim_{|z|\to0} \frac{F(x,t)}{|z|^{2-N}}
	=
	\frac{1}{N(N-2)\omega_N}.
\]
This completes the proof.
\Qed
\end{Proof}

\begin{pro}\label{pro:2uplow} \ 
Let $N=2$.
Suppose that $\xi$ is H\"older continuous
with exponent $\alpha  >1/2$.
Then for each $t\in(0,T)$ the function $F(x,t)$ satisfies
\[
	F(x,t) = 
	\frac{1}{2\pi} \log \Big( \frac{1}{|x-\xi(t)|} \Big)
		+ o \Big( \log \frac{1}{|x-\xi(t)|} \Big) 
	\quad \mbox{ as } x\to\xi(t).
\]
\end{pro}

\begin{Proof}
We fix $t\in(0,T)$ and set $z:=x-\xi(t)$.
Setting $N=2$ in (\ref{eq:change}), we have
\begin{equation}\label{eq:2change}
	\begin{aligned}
		F(x,t) 
		&=
		(4\pi)^{-1}
		\int_{|z|^2/4t}^\infty
			\sigma^{-1}
				\exp
				\left( - \left| 
					\sigma^{1/2} \frac{z}{|z|}
					+ \frac{1}{2}\frac{\xi(t)-\xi(s)}{(t-s)^{1/2}}
				\right|^2 \right) \,d\sigma
		\\&=:
		(4\pi)^{-1} I(z,t).
	\end{aligned}
\end{equation}
Here, we rewrite $I(z,t)$ as
\[
	I(z,t) =
	\int_{|z|^2/4t}^\infty 
			\sigma^{-1} e^{-\sigma}
				\exp
				\left( 
					- \sigma^{1/2} \frac{z}{|z|}
					\cdot \frac{\xi(t)-\xi(s)}{(t-s)^{1/2}} 
				\right)
				\exp
				\left( 
					- \frac{1}{4} \frac{|\xi(t)-\xi(s)|^2}{(t-s)} 
				\right)
	\, d\sigma.
\]
First, we claim that the function $F$ satisfies
\begin{equation}\label{eq:lsup}
	\limsup_{|z|\to0} 
	\frac{F(x,t)}{\log(1/|z|)}
	\leq \frac{1}{2\pi}.
\end{equation}
To show this, we give an upper bound of $I(z,t)$ as
\[
\begin{aligned}
	I(z,t)  
	& \leq
	\int_{|z|^2/4t}^\infty 
			\sigma^{-1} e^{-\sigma}
				\exp
				\left( 
					- \sigma^{1/2} \frac{z}{|z|}
					\cdot \frac{\xi(t)-\xi(s)}{(t-s)^{1/2}} 
				\right)
	\, d\sigma
	\\ & \leq
	\int_{|z|^2/4t}^\infty 
			\sigma^{-1} e^{-\sigma}
				\exp
				\left( 
					\frac{L}{4^{\alpha-(1/2)}}
					|z|^{2\alpha-1}\sigma^{1-\alpha}
				\right)
	\, d\sigma.
\end{aligned}
\]
For sufficiently small $|z|>0$,  we have
\[
\begin{aligned}
	I(z,t)  
	& \leq
	\int_1^\infty
	\sigma^{-1} e^{-\sigma+\sigma^{1-\alpha}}
	\, d\sigma
	+
	\exp
	\left(
		-\frac{|z|^2}{4t}
	\right)
	\exp
	\left(
		\frac{L}{4^{\alpha-(1/2)}}|z|^{2\alpha-1}
	\right)
	\int_{|z|^2/4t}^1 \sigma^{-1} 
	\, d\sigma
	\\ & =
	C(\alpha) 
	+
	\exp
	\left(
		-\frac{|z|^2}{4t}
	\right)
	\exp
	\left(
		\frac{L}{4^{\alpha-(1/2)}}|z|^{2\alpha-1}
	\right)
	\left(
		2\log\frac{1}{|z|} +\log(4t)
	\right)
\end{aligned}
\]
for some $C(\alpha)>0$.
Hence by (\ref{eq:2change}) and the above inequalities, we have
\[
\begin{aligned}
	\frac{F(x,t)}{\log(1/|z|)}
	& =
	\frac{I(z,t)}{4\pi \log(1/|z|)} 
	\\ & \leq
	\frac{1}{2\pi}
	\exp
	\left(
		-\frac{|z|^2}{4t}
	\right)
	\exp
	\left(
		\frac{L}{4^{\alpha-(1/2)}}|z|^{2\alpha-1}
	\right)
	\\ & \qquad +
	\frac
		{  C(\alpha) +
		\exp( -|z|^2/(4t) )
		\exp( 4^{(1/2)-\alpha} L |z|^{ 2\alpha-1} )
		\log(4t)  }
		{4\pi \log(1/|z|)}
	\\ & \to
	1/2\pi 
		\quad \mbox{ as } |z|\to 0.
\end{aligned}
\]
Consequently, we obtain (\ref{eq:lsup}).

Next, we claim that for any fixed $\eps\in(0,1)$ 
the function $F$ satisfies 
\begin{equation}\label{eq:linf}
	\liminf_{|z|\to0} 
	\frac{F(x,t)}{\log(1/|z|)}
	\geq \frac{1}{2\pi}(1-\eps).
\end{equation}
To show this, we give a lower bound of $I(z,t)$.
Now $(|z|^{2-\eps}/4t, |z|^\eps/4t) \subset (|z|^2/4t, \infty)$ holds.
Then, by using H\"older continuity, we directly calculate 
\[
\begin{aligned}
	I(z,t) 
	& \geq
	\int_{|z|^{2-\eps}/4t}^{|z|^\eps /4t}
			\sigma^{-1} e^{-\sigma}
				\exp
				\left( 
					- \sigma^{1/2} \frac{z}{|z|}
					\cdot \frac{\xi(t)-\xi(s)}{(t-s)^{1/2}} 
				\right)
				\exp
				\left( 
					- \frac{1}{4} \frac{|\xi(t)-\xi(s)|^2}{(t-s)} 
				\right)
	\, d\sigma
	\\ & \geq
	\int_{|z|^{2-\eps}/4t}^{|z|^\eps /4t}
		\sigma^{-1} e^{-\sigma}
			\exp
			\left( 
				- \frac{L}{4^{\alpha-(1/2)}}
				|z|^{2\alpha-1}
				\sigma^{1-\alpha}
			\right)
			\exp
			\left( 
				-\frac{L^2}{4^{2\alpha}}
				|z|^{4\alpha-2}
				\sigma^{-2\alpha+1}
			\right)
	\, d\sigma,
\end{aligned}
\]
where $L>0$ is a H\"older constant.
Since we assume $\alpha >1/2$, we have the following estimate:
\[
\begin{aligned}
	I(z,t)
	& \geq
	\exp
	\left(
		-\frac{|z|^\eps}{4t}
	\right)
	\exp
		\left( 
			- \frac{L}{4^{\alpha-\frac{1}{2}}}
			|z|^{2\alpha-1}
			\Big(
				\frac{|z|^\eps}{4t}
			\Big)^{1-\alpha}
		\right)
	\\ & \qquad \quad \times
	\exp
		\left( 
			-\frac{L^2}{4^{2\alpha}}
			|z|^{4\alpha-2}
			\Big(
				\frac{|z|^{2-\eps}}{4t}
			\Big)^{-2\alpha+1}
		\right)
	\int_{|z|^{2-\eps}/4t}^{|z|^\eps/4t}
		\sigma^{-1}
	\, d\sigma
	\\ & =
	\exp
	\left(
		-\frac{|z|^\eps}{4t}
	\right)
	\exp
		\left( 
			- \frac{L}{2} 
			t^{\alpha-1}
			|z|^{ 2\alpha-1 + \eps(1-\alpha) }
		\right)
	\\ & \qquad \quad \times
	\exp
		\left( 
			-\frac{L^2}{4}
			t^{2\alpha-1}
			|z|^{\eps(2\alpha-1)}
		\right)
	2(1-\eps) \log \frac{1}{|z|}.
\end{aligned}
\]
Hence by (\ref{eq:2change}) and  the above inequalities, we have
\[
\begin{aligned}
	\frac{F(x,t)}{\log(1/|z|)}
	& =
	\frac{I(z,t)}{4\pi \log(1/|z|)} 
	\\ & \geq
	\frac{1-\eps}{2\pi}
	\exp
	\left(
		-\frac{|z|^\eps}{4t}
	\right)
	\exp
		\left( 
			- \frac{L}{2}
			t^{\alpha-1}
			|z|^{ 2\alpha-1 + \eps(1-\alpha) }
		\right)
	\\ & \qquad \quad \times
	\exp
		\left( 
			-\frac{L^2}{4}
			t^{2\alpha-1}
			|z|^{\eps(2\alpha-1)}
		\right)
	\\ & \to
	(1-\eps)/2\pi 
		\quad \mbox{ as } |z|\to 0,
\end{aligned}
\]
so that (\ref{eq:linf}) holds.
These two claims imply that 
for any $\eps\in(0,1)$ the function $F$ satisfies
\[
	\frac{1-\eps}{2\pi}
	\leq
	\liminf_{|z|\to0} 
	\frac{F(x,t)}{\log(1/|z|)}
	\leq
	\limsup_{|z|\to0} 
	\frac{F(x,t)}{\log(1/|z|)}
	\leq 
	\frac{1}{2\pi}.
\]
Then
\[
	\lim_{|z|\to0} 
	\frac{F(x,t)}{\log(1/|z|)}
	=
	\frac{1}{2\pi}.
\]
This completes the proof.
\Qed
\end{Proof}

Now Theorem \ref{th:example} immediately follows from 
 Propositions \ref{pro:uplow}, \ref{pro:2uplow} 
and Proposition~\ref{pro:Fheat}.


\

\noindent
{\bf Acknowledgements.}
The authors would like to thank the referee for valuable comments.
The authors would also like to thank Dr. Toru Kan 
for many discussions. 
The second author was partially supported by 
JSPS KAKHENHI Grant-in-Aid for Challenging Exploratory Research (No.~22604020).


\providecommand{\bysame}{\leavevmode\hbox to3em{\hrulefill}\thinspace}
\providecommand{\MR}{\relax\ifhmode\unskip\space\fi MR }
\providecommand{\MRhref}[2]{%
  \href{http://www.ams.org/mathscinet-getitem?mr=#1}{#2}
}
\providecommand{\href}[2]{#2}

\end{document}